\documentclass[10pt]{article}

\usepackage[intlimits]{amsmath}

\usepackage{amssymb}

\usepackage[top=1.5in, bottom=1.5in, left=1.5in, right=1in]{geometry}

\newtheorem{theorem}{Theorem}[section]
\newtheorem{remark}{Remark}[section]
\newtheorem{definition}{Definition}[section]

\newtheorem{proposition}{Proposition}[section]
\newtheorem{lemma}{Lemma}[section]

\newtheorem{rule-def}[theorem]{Rule}
\newtheorem{example}{Example}[section]

\begin{document}

\title{\textbf{Existence uniqueness for a class of Nonlinear Discrete Boundary Value Problems}}
\author{Mandeep Singh$^a$, Amit K. Verma$^b$\thanks{Corresponding author Email: akviitp@gmail.com, $^a$mandeep04may@yahoo.in}
\\\small{\it{$^{a}$Department of Mathematics, DIT University, Dehradun, Uttarakhand-248009, India.}}
\\\small{\it{$^{b}$Department of Mathematics, IIT Patna, Bihta-801103, Bihar, India.}}
}
\date{}
\maketitle
\begin{abstract}
A monotone iterative method is proposed to solve nonlinear discrete boundary value problems with the support of upper and lower solutions. We establish some new existence results. Under some sufficient conditions, we establish maximum principle for linear discrete boundary value problem, which relies on Green's function and its constant sign. We then use it to establish existence of unique solution for the nonlinear discrete boundary value problem $\Delta^2 y(t-1)= f(t, y(t)),~t\in[1, T]$, $y(0)=0,~y(T+1)=0$.\\

\noindent\textit{Keywords:} Nonlinear Discrete BVPs; Monotone iterative technique; Upper and lower Solutions; Discrete Green's Function; Maximum Principle.\\

\noindent\textit{AMS Subject Classification:} 34B16; 39A10
\end{abstract}

\begin{abstract}
A monotone iterative method is proposed to solve nonlinear discrete boundary value problems with the support of upper and lower solutions. We establish some new existence results. Under some sufficient conditions, we establish maximum principle for linear discrete boundary value problem, which relies on Green's function and its constant sign. We then use it to establish existence of unique solution for the nonlinear discrete boundary value problem $\Delta^2 y(t-1)= f(t, y(t)),~t\in[1, T]$, $y(0)=0,~y(T+1)=0$.
\end{abstract}

\section{Introduction}
The main aspire of this work is to develop monotone iterative technique for the following nonlinear discrete boundary value problem
\begin{eqnarray}
\label{N-Discrete-P-1}&&-\Delta^2 y(t-1)= f(t, y(t)),\quad\quad  t\in[1, T],\\
\label{N-Discrete-P-2}&&y(0)=0,\quad~~ y(T+1)=0,
\end{eqnarray}
where $T$ is a positive integer, $[1, T]$ is the discrete interval $\{1,2,\cdots, T\}$, $f:[0,T]\times R\rightarrow R,$ and $\Delta$ is the forward difference operator. Here $f(t,y)$ is defined for all $t$ in $[1,T]$ and for all real numbers $y$. Also $f(t,y)$ satisfies growth conditions with respect to $y$ known as one sided Lipschitz condition given as
\begin{eqnarray*}
&& y\leq w \Longrightarrow  f(t,w)- f(t,y) \geq  M(w - y).
\end{eqnarray*}
We assume that $f(t,y)$ is continuous in $y$ for each $t\in[1,T]$.
Agarwal et al. \cite{Agarwal2004} used critical point theory and  discussed the existence result for the same nonlinear discrete boundary value problem (\ref{N-Discrete-P-1})-(\ref{N-Discrete-P-2}).

We introduce monotone iterative scheme for nonlinear discrete boundary value problem (\ref{N-Discrete-P-1})-(\ref{N-Discrete-P-2}) defined as,
\begin{multline}\label{Discrete-MIT}
-\Delta^2 y_{n+1}(t-1)-\lambda y_{n+1}(t)=f(t, y_{n}(t))-\lambda y_{n}(t),\\
\noindent y_{n+1}(0)=0,\quad y_{n+1}(T+1)=0, \lambda\in\mathbb{R} \setminus\{0\}.
\end{multline}
This technique is a discrete version of Picard type monotone iterative technique \cite{Picard1893}, (see \cite{Coster2006} and references their in). For discrete boundary value problems, we could observe only few references \cite{Wang1998, Wang1999, Zhuang1996, SR2011}.

We establish maximum principle for the corresponding linear discrete boundary value problem. We construct Green's function and prove that it is of constant sign for linear discrete boundary value problems (see \cite{Kelly2001}).

By using concept of upper and lower solutions with monotone iterative technique we establish a new existence result for nonlinear discrete boundary value problem. This existence result reveals that the upper and lower solutions are initial guess for monotone iterative technique, and generate monotonically convergent sequences. Thus we obtain the existence of unique solution for nonlinear discrete boundary value problem  (\ref{N-Discrete-P-1})-(\ref{N-Discrete-P-2}).
\section{Linear Discrete BVP}
To explore the existence results for nonlinear discrete boundary value problem (\ref{N-Discrete-P-1})-(\ref{N-Discrete-P-2}), we consider the following  non-homogeneous linear discrete boundary value problem
\begin{eqnarray}
\label{Discrete-NH-1}&&-Ly \equiv -\Delta^2 y(t-1)-\lambda y(t)=h(t),~~~~~~ t\in[1, T],\\
\label{Discrete-NH-2}&&y(0)=0,~~~~~~ y(T+1)=B,
\end{eqnarray}
where $B$ is any arbitrary constant. The corresponding homogeneous discrete boundary value problem will be
\begin{eqnarray}
\label{Discrete-H-1}&& \Delta^2 y(t-1)+\lambda y(t)=0,~~~~~~ t\in[1, T],\\
\label{Discrete-H-2}&&y(0)=0,~~~~~~ y(T+1)=0.
\end{eqnarray}
Solving non-homogeneous discrete boundary value problem (\ref{Discrete-NH-1})-(\ref{Discrete-NH-2}) is equivalently to find a $y$, such that
\begin{eqnarray}
y(t)= \bar{y} - \sum^{T}_{s=1}{G(t,s)h(s)}
\end{eqnarray}
where $\bar{y}$ is the solution of homogeneous difference equations (\ref{Discrete-H-1}), with non-homogeneous boundary conditions (\ref{Discrete-NH-2}) and $G(t,s)$ is the Green's function of (\ref{Discrete-H-1})-(\ref{Discrete-H-2}). Here, we discuss the solution of nonhomogeneous discrete boundary value problem (\ref{Discrete-NH-1})-(\ref{Discrete-NH-2}). We divide it into the following cases.
\begin{remark}\label{eigenvalue}
The characteristic equation for \eqref{Discrete-H-1} is
$$m^2+(\lambda-2)m+1=0.$$
If $|\lambda-2|<2,$ and $ \cos\theta = \frac{(2-\lambda)}{2},$ then $$m=\cos\theta\pm i \sin{\theta}=e^{\pm i\theta}.$$
Therefore a general solution of equation \eqref{Discrete-H-1} is
$$y(t)=c_1 \cos\left(\theta t\right)+c_2\sin\left(\theta t\right).$$
By using equation \eqref{Discrete-H-2}, we have
$$y(0)=c_1=0,\quad y(T+1)=c_2\sin\left({\theta}\right)=0,$$
which gives $$\theta=\frac{n\pi}{T+1}.$$
Thus all the eigenvalues are given by, $$\lambda_n=2-2\cos\left({\frac{n\pi}{T+1}}\right), \quad n=1,2,3\cdots T.$$
\end{remark}
\subsubsection{Case I: $|\lambda-2|<2$}
\begin{lemma}\label{Discrete-Lemma-1}
 The Green's function $G(t,s)$ of the discrete boundary value problem (\ref{Discrete-H-1})-(\ref{Discrete-H-2}), for $0<\lambda<\lambda_1$, is given by
\begin{eqnarray}\label{Discrete-GF-lambda-P}
G(t,s)=\left\{
  \begin{array}{ll}
 \frac{-\sin {\theta(T+1-s)}\sin{\theta}t}{\sin{\theta}\sin{\theta(T+1)}} , & 0 \leq t\leq s; \\ \vspace{0.2cm}\\

\frac{-\sin {\theta(T+1-s)}\sin{\theta}t}{\sin{\theta}\sin{\theta(T+1)}}+\frac{\sin{\theta(t-s)}}{\sin{\theta}} , & s \leq t\leq T+1;
  \end{array}
\right.
\end{eqnarray}
where $\theta=\tan^{-1}\left(\frac{\sqrt{4-(\lambda-2)^2}}{(2-\lambda)}\right)$, where $\lambda_1$ (see Remark \ref{eigenvalue}) is the first eigenvalue of $(\ref{Discrete-H-1})$-$(\ref{Discrete-H-2})$.
\end{lemma}
\textbf{Proof.} We define the Green's function as given below
\begin{eqnarray*}
G(x,t)=\left\{
  \begin{array}{ll}
    u(t,s) , & {0 \leq t\leq s\leq T+1;} \\
    v(t,s), & {0\leq s \leq t \leq T+1;}
      \end{array}
\right.
\end{eqnarray*}
where $u(\cdot, s)$ is the solution of the following discrete boundary value problem for each fixed $s\in[1,T]$
\begin{eqnarray}
\label{GF-1}&&Lu=0,\\
\label{GF-2}&&u(0,s)=0,\\
\label{GF-3}&&u(T+1, s)=- y(T+1,s),
\end{eqnarray}
and $v(t,s):=u(t,s)+y(t,s)$, where $y(\cdot,\cdot)$ is the Cauchy function for $Ly=0$, $\forall$ $t, s\in [1,T]$.
For each $s\in[1,T]$, $v(\cdot,s)$ is a solution of $Ly=0$ satisfying the boundary condition $y(T+1,s)=0.$

The Cauchy function for $Ly=0$ is given by
\begin{eqnarray}
y(t,s)= \frac{\sin {\theta(t-s)}}{\sin{\theta}}
\end{eqnarray}
where $\theta=\tan^{-1}\left(\frac{\sqrt{4-(\lambda-2)^2}}{(2-\lambda)}\right)$. Now from equations $(\ref{GF-1})$-$(\ref{GF-3})$, we get
\begin{eqnarray}
u(t,s)=\frac{-\sin {\theta(T+1-s)}\sin{\theta}t}{\sin{\theta}\sin{\theta(T+1)}},
\end{eqnarray}
for each fixed $s\in[1,T]$.

 Next
\begin{eqnarray}
v(t,s) =&&u(t,s)+y(t,s),\\
=&&\frac{-\sin {\theta(T+1-s)}\sin{\theta}t}{\sin{\theta}\sin{\theta(T+1)}}+\frac{\sin{\theta(t-s)}}{\sin{\theta}}.
\end{eqnarray}
This completes the construction of Green's function.

\begin{lemma}\label{Discrete-Lemma-2}
Let $y$ be a solution of non-homogeneous discrete boundary value problem (\ref{Discrete-NH-1})-(\ref{Discrete-NH-2}), then
\begin{eqnarray}
\label{NH-S-L-P}y(t)= \frac{B \sin{\theta t}}{\sin{\theta (T+1)}} - \sum^{T}_{s=1}{G(t,s)h(s)}.
\end{eqnarray}
\end{lemma}
\textbf{Proof.} Suppose $\bar{y}$ is the solution of $Ly=0,\quad t\in[1, T]$, subject to $y(0)=0,~~y(T+1)=B$,
and $G(t,s)$ is the solution of homogeneous discrete boundary value problem (\ref{Discrete-H-1})-(\ref{Discrete-H-2}). Then the discrete boundary value problem (\ref{Discrete-NH-1})-(\ref{Discrete-NH-2}) is analogous to
\begin{eqnarray*}
y(t)= \bar{y} - \sum^{T}_{s=1}{G(t,s)h(s)}.
\end{eqnarray*}
The general solution of difference equation $Ly=0$ is given by
\begin{eqnarray*}
\bar{y}(t)= c_1 \cos{\theta t}+c_2 \sin{\theta t},
\end{eqnarray*}
where $\theta=\tan^{-1}\left(\frac{\sqrt{4-(\lambda-2)^2}}{(2-\lambda)}\right)$.

Since
\begin{eqnarray*}
&&\bar{y}(0)=0,~{\mathrm{and}}~\bar{y}(T+1)=B,
\end{eqnarray*}
we get
\begin{eqnarray*}
&&c_1= 0,\\
&&c_2=\frac{B}{\sin{\theta(T+1)}}.
\end{eqnarray*}
Hence the discrete boundary value problem (\ref{Discrete-NH-1})-(\ref{Discrete-NH-2}) is equivalent to
\begin{eqnarray*}
y(t)= \frac{B \sin{\theta t}}{\sin{\theta (T+1)}} - \sum^{T}_{s=1}{G(t,s)h(s)}.
\end{eqnarray*}

Here we state four Lemmas \ref{Discrete-Lemma-3}, \ref{Discrete-Lemma-4}, \ref{Discrete-Lemma-1-lambda-0} and \ref{Discrete-Lemma-2-lambda-0} without proof. Proof of these are similar to the case I
\subsubsection{Case II: $\lambda<0$}
\begin{lemma}\label{Discrete-Lemma-3}
 The Green's function $G(t,s)$ of the discrete boundary value problem (\ref{Discrete-H-1})-(\ref{Discrete-H-2}) for $\lambda<0$, is given by
\begin{eqnarray}\label{Discrete-GF-lambda-N}
G(t,s)=\left\{
  \begin{array}{ll}
 \frac{1}{(\alpha-\beta)}\left({\alpha^{T+1-s}}-{\beta^{T+1-s}}\right)\frac{(\beta^t-\alpha^t)}{(\alpha^{T+1}-\beta^{T+1})} , & 0 \leq t\leq s; \\
 \vspace{0.2cm}\\
\frac{1}{(\alpha-\beta)}\left({\alpha^{T+1-s}}-{\beta^{T+1-s}}\right)\frac{(\beta^t-\alpha^t)}{(\alpha^{T+1}-\beta^{T+1})}+
\frac{(\alpha^{t-s}-\beta^{t-s})}{(\alpha-\beta)} , & s \leq t\leq T+1;
  \end{array}
\right.
\end{eqnarray}
\end{lemma}
 where $\alpha=\frac{(2-\lambda)+\sqrt{(\lambda-2)^2-4}}{2}$ and $\beta=\frac{(2-\lambda)-\sqrt{(\lambda-2)^2-4}}{2}$.

\begin{lemma}\label{Discrete-Lemma-4}
Let $y$ be the solution of non-homogeneous difference equation (\ref{Discrete-NH-1})-(\ref{Discrete-NH-2}), then
\begin{eqnarray}
\label{NH-S-L-N}y(t)= \frac{(\alpha^{t}-\beta^{t})B}{(\alpha^{T+1}-\beta^{T+1})} - \sum^{T}_{s=1}{G(t,s)h(s)}.
\end{eqnarray}
\end{lemma}
\subsubsection{Case III: $\lambda=0$}
\begin{lemma}\label{Discrete-Lemma-1-lambda-0}
 The Green's function $G(t,s)$ of the discrete boundary value problem (\ref{Discrete-H-1})-(\ref{Discrete-H-2}) for $\lambda=0$, is given by
\begin{eqnarray}\label{Discrete-GF-lambda-0}
G(t,s)=\left\{
  \begin{array}{ll}
 \frac{t(s-(T+1))}{T+1} , & 0 \leq t\leq s; \\
 \vspace{0.2cm}\\

\frac{s(t-(T+1))}{T+1}, & s \leq t\leq T+1.
  \end{array}
\right.
\end{eqnarray}
\end{lemma}
\begin{lemma}\label{Discrete-Lemma-2-lambda-0}
Let $y$ be the solution of non-homogeneous difference equation (\ref{Discrete-NH-1})-(\ref{Discrete-NH-2}), then
\begin{eqnarray}
\label{NH-S-L-0}y(t)= \frac{B}{{T+1}} - \sum^{T}_{s=1}{G(t,s)h(s)}.
\end{eqnarray}
\end{lemma}
\begin{remark}\label{Remark} Using Lemma \ref{Discrete-Lemma-2}, Lemma \ref{Discrete-Lemma-4} and Lemma \ref{Discrete-Lemma-2-lambda-0} the solution of non-homogeneous linear discrete boundary value problems  (\ref{Discrete-NH-1})-(\ref{Discrete-NH-2}) can be written as
\begin{eqnarray}
\label{NH-Sol}y(t)= B\psi(t) - \sum^{T}_{s=1}{G(t,s)h(s)}.
\end{eqnarray}
where $\psi(t)$ is defined as $\frac{\sin{\theta t}}{\sin{\theta (T+1)}}$ or $\frac{(\alpha^{t}-\beta^{t})}{(\alpha^{T+1}-\beta^{T+1})}$ or $\frac{B}{{T+1}}$ and $ G(t, s)$ is defined by \eqref{Discrete-GF-lambda-P} or \eqref{Discrete-GF-lambda-N} or \eqref{Discrete-GF-lambda-0}, respectively.
\end{remark}

\begin{remark}\label{Remark-1}
For each fixed value of $s\in[1,T]$, $u(\cdot, s)$ is the solution of (\ref{GF-1}) which satisfies the discrete boundary conditions $u(0,s)=0$ and u(T+1, s)= - y(T+1,s), where $y(\cdot,\cdot)$ is a Cauchy function and satisfies $y(s,s)=0$, $y(s+1,s)>0$. As difference equation (\ref{GF-1}) disconjugate on $[0,T+1]$, i.e,
$$u(t,s)<0,$$ for $t\in[1,T]$. Also $v(\cdot,s)$ is a solution of $Ly=0$ satisfying the boundary condition $v(T+1,s)=0$, and $$v(s,s)=u(s,s)+y(s,s)=u(s,s)<0,$$
then $$v(t,s)<0,$$ for $t\in[1,T]$. Therefore
$$G(t,s)<0,$$ for $t,s\in[1,T]$.

\end{remark}
\section{Maximum Principle}
\begin{proposition}\label{Maximum-Principle}
If $y$ satisfies the non-homogeneous linear discrete boundary value problem
\begin{eqnarray*}
&&-\Delta^2 y(t-1)-\lambda y(t)=h(t),~~~~~~ t\in[1, T],\\
&&y(0)=0,~~~~~~ y(T+1)=B,
\end{eqnarray*}
with $h(t)\geq0$ and $B\geq0$, then $y(t)\geq0$ for all $t\in[1,T]$ and $\lambda<\lambda_1.$
\end{proposition}
\textbf{Proof.} The proof is an immediate consequences of Remarks \ref{Remark}, \ref{Remark-1}.

\section{Nonlinear Discrete BVP}
In this section, we examine the existence results for nonlinear discrete boundary value problem, with the support of monotone iterative method and upper and lower solutions of the nonlinear discrete boundary value problem.

Let us first define the bounds of the solution of the nonlinear discrete boundary value problems
\begin{definition}
A function $\beta_0(t)$ is an upper solution of nonlinear discrete boundary value problem (\ref{N-Discrete-P-1})-(\ref{N-Discrete-P-2}) if it satisfies
\begin{gather}\label{Discrete-Upper-Sol.}\left.
\begin{aligned}
&-\Delta^2 \beta_0(t-1)\geq f(t, \beta_0(t)),\quad\quad  t\in[1, T],\\
&\beta_0(0)=0,\quad~~ \beta_0(T+1)\geq0.
\end{aligned}\right\}
\end{gather}
\end{definition}

\begin{definition}
A function $\alpha_0(t)$ is a lower solution of nonlinear discrete boundary value problem (\ref{N-Discrete-P-1})-(\ref{N-Discrete-P-2}) if it satisfies
\begin{gather}\label{Discrete-Lower-Sol.}\left.
\begin{aligned}
&-\Delta^2 \alpha_0(t-1)\leq f(t, \alpha_0(t)),\quad\quad  t\in[1, T],\\
&\alpha_0(0)=0,\quad~~ \alpha_0(T+1)\leq0.
\end{aligned}\right\}
\end{gather}
\end{definition}
\begin{theorem}\label{Discrete-Theorem-1}
If $f: D_0 \rightarrow R$ is continuous on $D_0:=\{(t,y) \in [0,T+1] \times R :\alpha_0\leq y \leq \beta_0\}$ in $y$ for each $t$ and there exist a constant $ M >0$ such that for all $(t,y),(t,w)\in D_0$
\begin{eqnarray}
\label{Discrete-Lipschitz-P}y\leq w \Longrightarrow  f(t,w)- f(t,y) \geq  M(w - y),
\end{eqnarray}
then the region $D_0$, contains at least one solution of the nonlinear discrete boundary value problem (\ref{N-Discrete-P-1})-(\ref{N-Discrete-P-2}). If a constant $\lambda\leq M$ is chosen such that $ \lambda < \lambda_1$ then the sequences $\{\beta_n\}$ generated by
\begin{eqnarray}
\label{Discrete-App-Scheme}&&-\Delta^2 y_{n+1}(t-1)-\lambda  y_{n+1}(t)=F(t,y_{n}(t)),~~~~y_{n+1}(0)=0,\quad y_{n+1}(T+1)=0,
\end{eqnarray}
where $F(x,y_{n}(t))=f(t,y_{n}(t))-\lambda y_{n}$, converges monotonically (non-increasing) and uniformly towards a solution ${{\widetilde{\beta}}(t)}$ of (\ref{N-Discrete-P-1})-(\ref{N-Discrete-P-2}). Similarly $\alpha_0$ leads to a non-decreasing sequences $\{\alpha_n\}$ converging to a solution ${\widetilde{\alpha}}(t)$. Any solution ${z}(t)$ in $D_0$ must satisfy
\begin{eqnarray*}
{\widetilde{\alpha}(t)\leq{z}(t) \leq\widetilde{\beta}}(t).
\end{eqnarray*}
\end{theorem}
\textbf{Proof.} Making the use of equations (\ref{Discrete-Upper-Sol.}) and (\ref{Discrete-App-Scheme})(for $n=0$)
\begin{gather}\label{Discrete-Existence-eq-1}\left.
\begin{aligned}
&-\Delta^2(\beta_0-\beta_1)(t-1)-\lambda(\beta_0-\beta_1)(t)\geq 0,\\
&(\beta_0-\beta_1)=0,\quad (\beta_0-\beta_1)(T+1)\geq 0.
\end{aligned}~~~~~~~~~~\right\}
\end{gather}
As $(\beta_0-\beta_1)$ satisfies the above equation (\ref{Discrete-Existence-eq-1}), with $h(t)\geq 0$, and $B\geq0$, then by the Proposition \ref{Maximum-Principle}, we have $\beta_0\geq \beta_1$.

As $M-\lambda\geq0$, using the equations (\ref{Discrete-Lipschitz-P}) and (\ref{Discrete-App-Scheme}), we have
\begin{eqnarray*}
-\Delta^2 \beta_{n+1}(t-1)\geq (M-\lambda)(\beta_n-\beta_{n+1})(t)+f(t, \beta_{n+1}(t)),
\end{eqnarray*}
and if $\beta_n-\beta_{n+1}\geq0$, then
\begin{eqnarray}
\label{Discrete-Existence-eq-2}&&-\Delta^2\beta_{n+1}(t-1) \geq f(t,\beta_{n+1}(t)).
\end{eqnarray}
Since $\beta_0 \geq \beta_1$, then by making the use of equations (\ref{Discrete-Existence-eq-2}) (for $n=0$) and (\ref{Discrete-App-Scheme}) (for $n=1$) we get
\begin{eqnarray*}
&&-\Delta^2(\beta_1-\beta_2)(t-1)-\lambda(\beta_1-\beta_2)(t)\geq 0,\\
&&~(\beta_1-\beta_2)(0)=0, ~~~~(\beta_1-\beta_2)(T+1)\geq0.
\end{eqnarray*}
In the view of Proposition \ref{Maximum-Principle}, we have $\beta_1\geq\beta_2$.

Now with the support of equations (\ref{Discrete-Lower-Sol.}) and (\ref{Discrete-App-Scheme}) (for $n=0$), we get
\begin{eqnarray*}
&&-\Delta^2(\beta_1-\alpha_0)(t-1)-\lambda(\beta_1-\alpha_0)(t)\geq 0,\\
&&~(\beta_1-\alpha_0)(0)=0, ~~~~(\beta_1-\alpha_0)(T+1)\geq0,
\end{eqnarray*}
which gives $\beta_1\geq\alpha_0$, (Proposition \ref{Maximum-Principle}).

To use mathematical induction, we assume that $\beta_{n+1}\leq\beta_{n},$ $\alpha_0\leq\beta_{n+1}$ and show that $\beta_{n+2}\leq\beta_{n+1}$ and $\alpha_0\leq\beta_{n+2}$ for all $n$. Now making the use of equations (\ref{Discrete-App-Scheme}) (for $n+1$) and (\ref{Discrete-Existence-eq-2})
\begin{eqnarray*}
&&-\Delta^2(\beta_{n+1}-\beta_{n+2})(t-1)-\lambda(\beta_{n+1}-\beta_{n+2})(t)\geq 0,\\
&&~(\beta_{n+1}-\beta_{n+2})(0)=0, ~~~~(\beta_{n+1}-\beta_{n+2})(T+1)\geq0,
\end{eqnarray*}
we have $\beta_{n+1}\leq\beta_{n}$ (Propositions \ref{Maximum-Principle}).

From equation (\ref{Discrete-App-Scheme}) (for $n+1$) and (\ref{Discrete-Lower-Sol.})
\begin{eqnarray*}
&&-\Delta^2(\beta_{n+2}-\alpha_0)(t-1)-\lambda(\beta_{n+2}-\alpha_0)(t)\geq 0,\\
&&~(\beta_{n+2}-\alpha_0)(0)=0, ~~~~(\beta_{n+2}-\alpha_0)(T+1)\geq0.
\end{eqnarray*}
Thus we have $\alpha_0\leq\beta_{n+2}$ (Proposition \ref{Maximum-Principle}) and hence
we have
\begin{eqnarray*}
\alpha_{0}\leq\ldots \leq\beta_{n+1}\leq \beta_{n}\leq \ldots \leq\beta_{2}\leq \beta_{1}\leq \beta_{0},
\end{eqnarray*}
and if we choose $\alpha_0$ as an initial iterate, then we easily get
\begin{eqnarray*}
\alpha_{0}\leq\alpha_{1}\leq\alpha_{2}\leq \ldots \leq\alpha_{n}\leq\alpha_{n+1} \leq \ldots \leq \beta_{0}.
\end{eqnarray*}
Finally we prove that $\alpha_{n}\leq \beta_{n}$ for all $n$. For this by assuming $\alpha_{n}\leq \beta_{n}$, we show that $\beta_{n+1}\geq \alpha_{n+1}$. From equation (\ref{Discrete-App-Scheme})  it is easy to get
\begin{eqnarray*}
&&-\Delta^2(\beta_{n+1}-\alpha_{n+1})(t-1)-\lambda(\beta_{n+1}-\alpha_{n+1})(t)\geq 0,\\
&&~~(\beta_{n+1}-\alpha_{n+1})(0)=0,~~~~(\beta_{n+1}-\alpha_{n+1})(T+1)\geq0.
\end{eqnarray*}
Hence by Proposition \ref{Maximum-Principle}, $\beta_{n+1}\geq \alpha_{n+1}$. Thus we have
\begin{eqnarray*}
\alpha_{0}\leq\alpha_{1}\leq\alpha_{2}\leq \ldots \leq\alpha_{n}\leq\alpha_{n+1} \leq \ldots \leq\beta_{n+1}\leq \beta_{n}\leq \ldots \leq\beta_{2}\leq \beta_{1}\leq \beta_{0}.
\end{eqnarray*}
So the sequences $\beta_{n}$ and $\alpha_{n}$ are monotonically non-increasing and non-decreasing, respectively and are bounded by $\beta_0$ and $\alpha_0$. Hence by Dini's theorem they converges uniformly. Let $\beta(t)=\displaystyle\lim_{n\to\infty}\beta_{n}(t) $ and $\alpha(t)=\displaystyle\lim_{n\to\infty}\alpha_{n}(t)$.

The solution $\beta_{n+1}$ of equation (\ref{Discrete-App-Scheme}) is given by (Remark \ref{Remark-1}).
\begin{eqnarray*}
\beta_{n+1}= B\psi(t)  - \sum^{T}_{s=1}{G(t,s)(f(t,\beta_n(t))-\lambda \beta_n)}.
\end{eqnarray*}
Now by dominated convergence theorem, as $n$ approaches to $\infty$, we get
\begin{eqnarray*}
&&{{\widetilde{\beta}}(t)}= B\psi(t) - \sum^{T}_{s=1}{G(t,s)(f(t,{{\widetilde{\beta}}(t)})-\lambda {{\widetilde{\beta}}(t)})}.
\end{eqnarray*}
Which is the solution of boundary value problem (\ref{N-Discrete-P-1})-(\ref{N-Discrete-P-2}).

It is clear that any arbitrary solution $z(t)$ can be treated as upper solution $\beta_0(t)$, i.e., we get $z(t)\geq \alpha_0(t)$, similarly one concludes that $z(t)\leq \beta_0(t)$.
\begin{theorem}\label{Uniquness},
Let $f(t, y)$ is continuous in $y$ for each $t$ in $[1, T]$ and there is a constant $M>0$ such that
\begin{eqnarray}
\label{Lip_Unique}f(t,w)- f(t,y) \geq  M(w - y),
\end{eqnarray}
and $M<\lambda_1$. Then the nonlinear discrete boundary value problem \eqref{N-Discrete-P-1}--\eqref{N-Discrete-P-2} has unique solution.
\end{theorem}
\textbf{Proof.} The Proof of this theorem is an immediate consequence of Maximum principle (Proposition \ref{Maximum-Principle}) and both side Lipschitz condition \eqref{Lip_Unique}.
\section{Numerical Illustration}
\begin{example}
Consider the nonlinear discrete  boundary value problems
\begin{eqnarray}
\label{Discrete-Exp-1-eqn-1}&&-\Delta^2 y(t-1)= \frac{e^{y(t)}}{e^{(T+1)^2}}, ~~~t\in[1, T],\\
\label{Discrete-Exp-1-eqn-2}&& y(0)=0,~~y(T+1)=0.
\end{eqnarray}
\end{example}
Here, $\alpha_0=0$ and $\beta_0= (T+1)t-\frac{t^2}{2}$ are defined as lower and upper solutions of the solution of nonlinear discrete boundary value problem (\ref{Discrete-Exp-1-eqn-1})-(\ref{Discrete-Exp-1-eqn-2}), respectively. The nonlinear source term is continuous for all values of $y(t)$ and satisfies one sided Lipschitz condition, with constant $M=\frac{1}{e^{(T+1)^2}}$. By Theorem \ref{Discrete-Theorem-1} and Theorem \ref{Uniquness}, \eqref{Discrete-Exp-1-eqn-1}--\eqref{Discrete-Exp-1-eqn-2} has a unique solution.

\begin{example}
Consider the nonlinear discrete  boundary value problems
\begin{eqnarray}
\label{Discrete-Exp-2-eqn-1}&&-\Delta^2 y(t-1)= e^{t}-e^{y(t)}, ~~~t\in[1, T],\\
\label{Discrete-Exp-2-eqn-2}&& y(0)=0,~~y(T+1)=0.
\end{eqnarray}
\end{example}
Here, $\alpha_0=0$ and $\beta_0= t$ are defined as lower and upper solutions of the solution of nonlinear discrete boundary value problem (\ref{Discrete-Exp-2-eqn-1})-(\ref{Discrete-Exp-2-eqn-2}), respectively. The nonlinear source term is continuous for all values of $y(t)$ and satisfies one sided Lipschitz condition, with constant $M=e^{(T+1)}$. By Theorem \ref{Discrete-Theorem-1}, and Theorem \ref{Uniquness}, \eqref{Discrete-Exp-2-eqn-1}--\eqref{Discrete-Exp-2-eqn-2} has a unique solution.

\bibliography{msakv}
\end{document}